\documentclass[preprint,11pt]{elsarticle}
\usepackage{mathrsfs}
\usepackage{graphics}
\usepackage{latexsym,amssymb,amsmath,amscd,amscd,amsthm,amsxtra}

\newtheorem{thm}{Theorem}
\newtheorem{lem}{Lemma}
\newtheorem{cor}{Corollary}

\newtheorem{prop}{Proposition}

\begin{document}
\begin{frontmatter}

\title{On jumped Wenger graphs}

\author[wlp]{Li-Ping Wang}

\ead{wangliping@iie.ac.cn}

\author[wan]{Daqing Wan}
\ead{dwan@math.uci.edu}

\author[wangw]{Weiqiong Wang}
\ead{wqwang@chd.edu.cn}

\author[zhou]{Haiyan Zhou \corref{cor1}}
\ead{haiyanxiaodong@gmail.com}

\cortext[cor1]{Corresponding author}  

\address[wlp]{State Key Laboratory of Information Security, Institute of Information Engineering, \\
Data Assurance and Communications Security Research Center,\\
Chinese Academy of Sciences, Beijing 100093,   China}

\address[wan]{Department of Mathematics, University of California, Irvine, CA 92697-3875, USA}

\address[wangw]{School of Science, Chang'an University, Xi'an 710064, China}

\address[zhou]{School of Mathematics, Nanjing Normal University, Nanjing 210023, China}

\begin{abstract}
In this paper  we introduce a new infinite class of bipartite graphs, called jumped Wenger graphs, which are closely related to  Wenger graphs.    An tight upper bound of the diameter  and  the exact  girth of a jumped Wenger graph $J_m(q, i, j )$ for  integers $i, j$, $1\leq i <j \leq m+2$, are determined. In particular,  the exact diameter of the jumped Wenger graph $J_m(q, i, j)$ if $(i, j)=(m,m+2), (m+1,m+2)$ or $(m,m+1)$ is also obtained.    
\end{abstract} 



\end{frontmatter}

\section{Introduction} 

Recently,  researchers focus on a class of bipartite graphs related to Wenger graphs because of their nice graphical theory properties \cite{cll,clw3,lu2, lw, SHS,Viglione0,Viglione, Wenger}.
First we introduce these graphs.

Let $\mathbb{F}_q$ be a finite field of order $q=p^e$ where  $p$ is a prime and $e$ is a positive integer.  Let $m$ be a positive integer and let  $\mathfrak{P}=\mathbb{F}_{q}^{m+1}$ and  $\mathfrak{L}=\mathbb{F}_{q}^{m+1}$ be two copies of the $(m+1)$-dimensional vector space over  $\mathbb{F}_q$, which are called the point set and the line set respectively.  
Let $\mathfrak{G} =(V,E)$ be the bipartite graph with the vertex set $V=\mathfrak{P}\cup \mathfrak{L}$ and the edge set $E$ is defined as follow: there is an edge from a point $P=(p_1,p_2,\cdots,p_{m+1})\in \mathfrak{P}$  to a line $L=[l_1,l_2,\cdots,l_{m+1}]\in \mathfrak{L}$, denoted by $PL$,  if the following $m$ equalities hold:
\begin{eqnarray}\label{defnGraph}
  l_2+p_2 &=& l_1\,f_2(p_{1})  \nonumber  \\
 l_3+p_3 &=& l_1\,f_3(p_{1}) \nonumber  \\
 \vdots&\vdots& \vdots\\
 l_{m+1}+p_{m+1} &= & l_1\,f_{m+1}(p_{1}), \nonumber
\end{eqnarray} 
where $f_i(x)\in \mathbb{F}_q[x]$ for $2\leq i\leq m+1$.

If $(1,f_2(x),\ldots,f_{m+1}(x))=(1,x,x^2,\ldots,x^{m})$,  this class of graphs is called 
Wenger graphs \cite{cll}.  If  $(1,f_2(x),\ldots,f_{m+1}(x))=(1,x,x^{p},\ldots,x^{p^{m-1}})$, this class of graphs is called linearized Wenger graphs \cite{clw3}. In this paper, we focus on the more general case when  
$(1,f_2(x),\ldots$, $f_{m+1}(x))=(1,x,\ldots$, $ x^{i-1}, x^{i+1},\ldots$, $ x^{j-1},x^{j+1}, \ldots, x^{m+2})$, with $1\leq i<j\leq m+2$, which is called the jumped Wenger graphs with jump points at $x^{i}$ and $x^{j}$, denoted by $J_m(q, i, j )$.   In particular,  
if $j=m+2$, that is,  $(1,f_2(x),\ldots$, $f_{m+1}(x))=(1,x,\ldots, x^{i-1}, x^{i+1}$, $\ldots, x^{m}, x^{m+1})$, the graphs has only one jump point at $x^{i}$.   If $i=m+1$ and $j=m+2$, the graphs become Wenger graphs, denoted by $W_m(q)$.  

It is easy to get some properties from \cite{clw3}.  

 {\prop The jumped Wenger graph $J_{m}(q, i, j)$ for  $1\leq i<j \leq m+1$ is $q$-regular. }

{
\prop  \label{pro:connect}
If $m+2<q$, the jumped Wenger graph $J_{m}(q, i, j)$ for  $1\leq i<j \leq m+1$ is connect. 

}

In Section 2  we  give an upper bound of  the diameter of a jumped Wenger graph $J_m(q,i, j)$ for any integers $i, j$, $1\leq i <j\leq m+2$  and  the exact diameter for some special jumped Wenger graphs such as $(i, j)=(m,m+2), (m+1,m+2)$ or $(m,m+1)$. In Section 3, we determine the girth of  a jumped Wenger graph $J_m(q,i, j)$ for $1\leq i <j\leq m+2$.  Finally, we give the conclusions. 
  
\section{The diameter of jumped Wenger graphs }

Recall that a sequence of distinct vertices $v_1,\cdots,v_s$ in a simple connected graph  $\mathfrak{G}=(V,E)$ defines a {\it path} of length $s-1$ if $(v_i,v_{i+1})\in E$ for every $i, 1\leq i\leq s-1$. The {\em distance} between $v_i$ and $v_j$ is the number of edges in a shortest path joining $v_i$ and $v_j$. The {\em diameter} of a connected graph $\mathfrak{G}$ is the maximum distance between any two vertices of $\mathfrak{G}$. In  \cite{Viglione} it is shown that the diameter of the Wenger graph $W_m(q)$ is $2m+2$ if $1\leq m\leq q-1$.  The result for the linearized Wenger graphs also hold  if $m\leq e$ in \cite{cll}. By Proposition \ref{pro:connect},  for any $1\leq i< j\leq m+2$, the jumped Wenger graph $J_m(q, i, j)$  is connected if $m+2<q$.  So  we consider  the diameter of  these jumped Wenger graphs  in this section.       
 
Assume that $L=[l_1,\ldots,l_{m+1}]$ and $L'=[l'_1,\ldots, l'_{m+1}]$ are two lines and they  share a point 
$P=(p_1,\ldots,p_{m+1})$. Then we have 
\begin{eqnarray*}
l_k+p_k & = & l_1 f_{k}(p_{1}), k=2,\ldots,m+1 \\
l'_k+p_k & = & l'_1 f_{k}(p_{1}), k=2,\ldots,m+1,
\end{eqnarray*}
\noindent and so we get 
\begin{equation*} \label{eq:line}
l_k-l'_k  = (l_1-l'_1)f_{k}(p_1), k=2,\ldots,m+1. 
\end{equation*}

We obtain the following fact which will be used later. 

\medskip 

\noindent \textbf{Fact 1}.  $L$ and $L'$ with a common point are the same line if $l_1=l'_1$.

Assume $P=(p_1,\ldots,p_{m+1})$ and $P'=(p'_1,\ldots,p'_{m+1})$ be two points, both of which share a line 
$L=[l_1,\ldots,l_{m+1}]$.  Then we have 
\begin{eqnarray*}
l_k+p_k & = & l_1 f_{k}(p_{1}), k=2,\ldots,m+1\\
l_k+p'_k &= & l_1 f_{k}(p'_1), k=2,\ldots,m+1,
\end{eqnarray*} 
and so we get  
\begin{equation*}
p_k-p'_k =  l_1 (f_{k}(p_1)-f_{k}(p'_1)), k=2,\ldots,m+1.
\end{equation*}

In this case we have  $P-P'=(p_1-p'_1, l_1(p_1-p'_1), \ldots,  l_1(p^{i-1}_1-(p'_1)^{i-1})$, $ l_1(p^{i+1}_1-(p'_1)^{i+1}),\ldots,  l_1 (p^{j-1}_1-(p'_1)^{j-1}),\ldots, l_1 (p^{m+2}_1-(p'_1)^{m+2}))$ and so we obtain the following fact.  
\medskip

\noindent \textbf{Fact 2}. The  $P$ and $P'$ sharing the same line are  identical if $p_1=p'_1$.  

 \medskip 
Next, we recall elementary symmetric polynomials.   
If $f(x)=\prod_{i=1}^{n}(x-x_i)=x^n-a_1x^{n-1}+a_2x^{n-2}+\ldots+(-1)^{n}a_n$, then the $k$th 
  elementary symmetric polynomials 
  $ \sigma_{k}(x_1,\ldots,x_n)=\sum_{1\leq i_1<\ldots <i_{k}\leq n} x_{i_1}\cdots x_{i_k}$, $k=1,\ldots,n$.  In addition,  we define $\sigma_0(x_1,\ldots,x_n)=1$ and $\sigma_{-1}(x_1,\ldots,x_n)=0$.  
 So  it is straightforward to have the following iterative identity  
 \begin{equation} \label{eq:newton}
 \sigma_k(x_1,\ldots,x_n)=x_{n}\sigma_{k-1}(x_1,\ldots,x_{n-1})+\sigma_k(x_1,\ldots,x_{n-1}).
 \end{equation} 
  
For simplicity, we denote      
$$\sigma_{i,j}(x_1,\ldots,x_n)=\sigma_{i}(x_1,\ldots,x_n)\sigma_{j+1}(x_1,\ldots,x_n)-\sigma_{j}(x_1,\ldots,x_n)\sigma_{i+1}(x_1,\ldots,x_n), $$  
$M_{l,i,j}(x_1,\ldots,x_n)=\left(\begin{array}{cccc}
             1  & 1 & \ldots  & 1 \\
             x_1 & x_2 & \ldots & x_{n}\\
             \multicolumn{4}{c}\dotfill\\
             x^{i-1}_1 & x^{i-1}_2 & \ldots & x^{i-1}_n \\
             x^{i+1}_1 & x^{i+1}_2 & \ldots & x^{i+1}_n \\
              \multicolumn{4}{c}\dotfill\\
              x^{j-1}_1 & x^{j-1}_2 & \ldots & x^{j-1}_n\\
             x^{j+1}_1 & x^{j+1}_2 & \ldots & x^{j+1}_n \\
             \multicolumn{4}{c}\dotfill\\
             x^{l+1}_1 & x^{l+1}_2 & \ldots & x^{l+1}_n
           \end{array}\right )$ and \\
            $V_{l}(x_1,\ldots,x_n)=\left(\begin{array}{cccc}
             1  & 1 & \ldots  & 1 \\
             x_1 & x_2 & \ldots & x_{n}\\
             {x}^{2}_1 & {x}^{2}_2 & \ldots & {x}^{2}_n \\
             \multicolumn{4}{c}\dotfill\\
             x^{l-1}_1 & x^{l-1}_2 & \ldots & x^{l-1}_n
           \end{array}\right )$ for a vector $(x_1,\ldots,x_n)\in \mathbb{F}^{n}_{q}$ and $0\leq i<j\leq l+1$.

 {\lem \label{lem-del2} For a vector  $(x_1,\ldots, x_{n})\in \mathbb{F}^n_q$ and $1\leq i<j\leq n+1$, we have 
$$|M_{n,i,j}(x_1,\ldots,x_n)|=(-1)^{i+j-1}\sigma_{n-i,n-j}(x_1,\ldots,x_n)
           \prod_{1\leq k<l \leq n}(x_{l}-x_{k})$$
 and so the determinant is not $0$ if and only if the  $x_i$ are distinct and  
  $\sigma_{n-i,n-j}(x_1$, $\ldots$, $x_n)\neq 0$. }
 
\noindent \textbf{Proof.}  Consider the following determinant 
\begin{eqnarray*} 
|V_{n+2}(x_1,\ldots,x_n,x,y)|
  & = &  (y-x)\prod_{k=1}^{n}(y-x_k)\prod_{k=1}^{n}(x-x_k)\prod_{1\leq k<l\leq n}(x_l-x_k)\\   
&  = &  (y-x)\sum_{k=0}^{n}(-1)^k\sigma_k(x_1,\ldots,x_n)y^{n-k} \\
& & \times  \sum_{k=0}^{n}(-1)^k\sigma_k(x_1,\ldots,x_n)x^{n-k} \prod_{1\leq k<l\leq n}(x_l-x_k).  
 \end{eqnarray*}
  
\noindent  So $(-1)^{i+j+1}|M_{n,i,j}(x_1,\ldots,x_n)|$ is the coefficient of $x^iy^j$ in $|V_{n+2}(x_1,\ldots, x_n$, $x,y)|$, i.e., $(-1)^{i+j+1}\sigma_{n-i,n-j}(x_1,\ldots,x_n)\prod_{1\leq k<l \leq n}(x_{l}-x_{k})$.  \hfill $\Box$
  
 \medskip 
 
  If $j=n+1$, it is reduced to the case of just deleting $x^i$ and so we get the following corollary.     

{\cor
For any $i$ with $1\leq i\leq n$, we have 
\[|M_{n,i,n+1}(x_1,\ldots,x_n)|= (-1)^{i+n}\sigma_{n-i}(x_1,\ldots,x_n)\prod_{1\leq k<l \leq n}(x_{l}-x_{k}).\] }

 {\lem \label{lem-2.2}
 Let $q>n+1$.  For any integer $i$, $0\leq i\leq n$, there are distinct $x_1,\ldots,x_n\in \mathbb{F}_q$ 
 such that $\sigma_{n-i}(x_1,\ldots,x_n)\neq 0$. 
 }
 
\medskip 

\noindent \textbf{Proof.}  We proceed by induction on $n-i$. 
If $n-i=0$, that is, $i=n$, the result holds obviously since $\sigma_0(x_1,\ldots,x_n)=1$. 
Suppose that the result holds for $n-i-1$. 
Since $\sigma_{n-i}(x_1,\ldots,x_n)=x_n\sigma_{n-i-1}(x_1,\ldots$, $x_{n-1})+\sigma_{n-i}(x_1,\ldots,x_{n-1})$,  by induction hypothesis, there are distinct $x_1,\ldots,x_{n-1}$ such that  $\sigma_{n-i-1}(x_1,\ldots,x_{n-1})\neq 0$.  Because of  $q>n+1$, one can choose  $x_n$ different from all other $x_i$, $1\leq i\leq n-1$,  such that $\sigma_{n-i}(x_1,\ldots,x_n)\neq 0$.  \hfill $\Box$

{\lem \label{lem-2.3} 
 Let $q>n+1$.  For any integers $i, j$,  $0\leq i< j\leq n+1$, there are distinct $x_1,\ldots,x_n\in \mathbb{F}_q$ such that  $\sigma_{n-i, n-j}(x_1,\ldots,x_n)\neq 0$.   }
  
\medskip  

 \noindent \textbf{Proof.} By Lemma \ref{lem-2.2}, the result holds for only induction on $n-i$ or $n-j$.  Then  we proceed by induction on  both $n-i$ and $n-j$.   If  $n-i=0$ and $n-j=-1$, we have  $\sigma_{0,-1}(x_1,\ldots,x_n)=\sigma_{0}(x_1,\ldots,x_n)\sigma_{0}(x_1,\ldots,x_n)-\sigma_{1}(x_1,\ldots,x_n)\sigma_{-1}(x_1,\ldots,x_n)$ $=1$ and so the result holds . Suppose that the result holds for any $n-j-1$ and $n-i-1$.  By Eq. (\ref{eq:newton}), we have 
 \begin{eqnarray}\label{eq:recursion}
 \sigma_{n-i, n-j}(x_1,\ldots,x_n)&= & x^2_{n}\sigma_{n-i-1,n-j-1}(x_1,\ldots,x_{n-1})+x_n(\sigma_{n-i-1}(x_1,\ldots,  \nonumber \\
 & & x_{n-1})\sigma_{n-j+1}(x_1, \ldots,x_{n-1}) -\sigma_{n-i+1}(x_1,\ldots,x_{n-1}) \nonumber \\
 & & \sigma_{n-j-1}(x_1,\ldots,x_{n-1}))+\sigma_{n-i,n-j}(x_1,\ldots,x_{n-1}). 
 \end{eqnarray} 
By induction hypothesis, there are distinct $x_1,\ldots,x_{n-1}\in \mathbb{F}_q$ such that 
\[\sigma_{n-i-1,n-j-1}(x_1, \ldots, x_{n-1})\neq 0\]
\noindent  and so one can choose $x_n$  not only different from  $x_1$,
$\ldots, x_{n-1}$  but also  different from any roots of the above quadratic 
Eq. (\ref{eq:recursion}) since $q>n+1$.    \hfill $\Box$

{\lem \label{lem-2.4}
 Let $q>n+1$.  For any fixed $x_1\in \mathbb{F}_q$, any integer $0\leq i< j\leq n+1$, there are distinct $x_2,\ldots,x_n\in \mathbb{F}_q$, which are also different from $x_1$, such that 
 $\sigma_{n-i, n-j}(x_1,\ldots,x_n)\neq 0$. }
  
\medskip 

 \noindent \textbf{Proof.}   We also proceed by induction on both  $n-i$ and $n-j$.    If  $n-i=0$ and $n-j=-1$, we have  $\sigma_{0,-1}(x_1,\ldots,x_n)=1$ and so the result holds clearly. 
 Suppose the result holds for  $n-j-1$ and $n-i-1$ for fixed $x_1$.  By Eq. (\ref{eq:newton}), we have 
 \begin{eqnarray} \label{eq:fix} 
 \sigma_{n-i, n-j}(x_1,\ldots,x_n)&= & x^2_n\sigma_{n-i-1, n-j-1}(x_1,\ldots,x_{n-1})+x_n(\sigma_{n-i-1}(x_1,\ldots, \nonumber \\
 & & x_{n-1})\sigma_{n-j+1}(x_1, \ldots,x_{n-1}) -\sigma_{n-i+1}(x_1, \ldots,x_{n-1}) \nonumber\\
 & & \sigma_{n-j-1}(x_1,\ldots,x_{n-1}))+\sigma_{n-i, n-j}(x_1,\ldots,x_{n-1}). 
 \end{eqnarray} 
By induction hypothesis, there are distinct $x_2,\ldots,x_{n-1}$, which is also different from $x_1$, such that $\sigma_{n-i-1, n-j-1}(x_1,\ldots$, $x_{n-1})\neq 0$ and so one can choose $x_n$ not only different from  $x_1$, $x_2$, 
$\ldots, x_{n-1}$  but also  different from any roots of the above quadratic 
Eq. (\ref{eq:fix}) since $q>n+1$.  \hfill $\Box$

{\thm \label{thm-5} If $1\leq m< q-2$, the diameter  of the jumped Wenger graph $J_m(q, i, j)$ for any integers $i,j$,  $1\leq i<j \leq m+1$,  is at most $2(m+1)$. }

\medskip

\noindent \textbf {Proof.} First,  we consider the distance between any two vertices $L$ and $L'$ in $\mathfrak{L}$  of the jumped Wenger graph $J_m(q, i, j)$. If $L_1P_1\ldots P_sL_{s+1}$ is a path in $J_m(q, i, j)$ between  $L=L_1$ and $L'=L_{s+1}$, where $L_h=[l_1^{(h)},\cdots,l_{m+1}^{(h)}]$ and
$P_h=(p_1^{(h)},\cdots,p_{m+1}^{(h)})$, we have
\begin{eqnarray*} \label{eq:path}
  l_k^{(h+1)}-l_k^{(h)} & =& (l_1^{(h+1)}-l_1^{(h)})f_{k}(p^{(h)}_1), k=2,\cdots,m, h=1,\cdots,s.\\
   l_{m+1}^{(h+1)}-l_{m+1}^{(h)} &= & (l_1^{(m+1)}-l_1^{(m)})f_{m+1}(p_1^{(h)}), h=1,\cdots, s. 
\end{eqnarray*}
 Therefore  there are elements $t_h =l_1^{(h+1)}-l_1^{(h)} $, $x_h = p_1^{(h)} \in \mathbb{F}_q$, $1\leq h\leq s$,  such that 
\begin{equation}\label{f-11}
  (L_{s+1}-L_1)^T =M_{m+1,i,j}(x_1,\ldots,x_s)(t_1\ \ t_2\ \ \ldots \ \ t_s )^T.
\end{equation} 
\noindent In addition, by  Fact 1 and 2, $x_1,\ldots,x_{s}$ should be distinct and $t_i\neq 0$ for all $i$, $1\leq i\leq s$.

For any integers  $i, j$, $1\leq i< j\leq m+2$,  
take $s=m+1$ in Eq. (\ref{eq:path}).  By Lemma \ref{lem-del2} and \ref{lem-2.3}, there are  distinct $x_1, \ldots, x_{m+1} \in \mathbb{F}_q$ such that  $|M_{m+1, i, j}(x_1,\ldots,x_{m+1})|\neq 0$, and  thus Eq. (\ref{f-11}) has a unique solution for $t_1, t_2, \ldots, t_s$.  If some $t_i=0$, then $P_i$ and $L_{i+1}$ can be deleted from the path.  Thus the distance of any two vertices in $\mathfrak{L}$  is at most  $2(m+1)$. 

 Then,  consider any vertices $P$ and $P'$ in  $\mathfrak{P}$.  Let $P_1L_1\ldots L_sP_{s+1}$ be a path in  $J_{m}(q, i, j)$ between $P=P_1$ and $P'=P_{s+1}$, where $L_h=[l^{(h)}_1,\ldots, l^{(h)}_{m+1}]$ and 
 $P_h=(p^{(h)}_1,\ldots, p^{(h)}_{m+1})$.  Therefore, we have 
  \begin{equation}\label{f-12}
  (P_{s+1}-P_1)^T =M_{m+1, i, j}(x_1,\ldots,x_s)(t_1\ \ t_2\ \ \ldots \ \ t_s )^T,
\end{equation}
 where $x_h = p_1^{(h)}$, $1\leq h\leq s$, $t_h=l_1^{(h)}-l_1^{(h+1)}$, $h=2,\ldots,s-1$, $t_1=-l^{(1)}_1$, 
 $t_{s}=l^{(s)}_1$.  Similarly,  take $s=m+1 $. By Lemma \ref{lem-del2} and Lemma \ref{lem-2.3}, there are  distinct $x_1, \ldots, x_{m+1} \in \mathbb{F}_q$ such that  $|M_{m+1, i, j}(x_1,\ldots,x_{m+1})|\neq 0$, and  thus Eq. (\ref{f-12}) has a unique solution for $t_1, t_2, \ldots, t_s$.  Thus the distance of any two vertices in $\mathfrak{P}$  is at most  $2(m+1)$. 
 
    Next, we consider the distance between a vertex $P=(p_1,\ldots,p_{m+1})\in \mathfrak{P}$ and a vertex $L \in \mathfrak{L}$.  First we choose a line $L_1\in \mathfrak{L}$ such that it is adjacent to $P$.  From the earlier discussion,  there exists a path  from $L_1$ to $L=L_{s+1}$ with distance at most $2(m+1)$. We modify the earlier construction so that the path goes through the vertex $P$. Namely,  In Eq. (\ref{f-11}), we let
 $x_1=p_1$ and choose distinct  $x_2,\ldots,x_n$ so that  $|M_{m+1, i, j}(x_1,\ldots,x_{m+1})|\neq 0$ by Lemma \ref{lem-del2} and \ref{lem-2.4} . Then there is a unique solution $t_1,\ldots,t_{m+1}$ and so there is a path between $L_1$ and $L$ with length at most $2(m+1)$ passing through $P$ by Fact 1. Therefore the distance of $P$ and $L$ is less than  or  equal to $2(m+1)$.

  Finally, we obtain the same result for the distance between a line $L$ and a point $P$.    \hfill $\Box$

{\thm   If $(i, j)= (m,m+2), (m+1,m+2), (m,m+1)$, then the diameter of $J_m(q, i, j)$ is $2(m+1)$.  } 

\medskip 

\noindent \textbf {Proof.}  By Theorem 1, it suffices to show that there  are $L$ and $L'$ in $\mathfrak{L}$ with  $ L-L'=[0,\ldots,0,1]$ such that  the distance  between them is at least $2(m+1)$.  
Otherwise,  there is  a shortest path $LP_1\ldots P_s L'$ between $L$ and $L'$ with $s\leq m$, where $P_h=(p^{(h)}_1,\ldots, p^{(h)}_{m+1})$ and $x_h=p^{(h)}_1$ for $1\leq h\leq s$.  The Eq. (\ref{f-11}) is easily reduced to 
\begin{equation} \label{eq:reduce}
(L'- L)^T=M_{s', i, j}(x'_1,\ldots, x'_{s'})(t_1,\ldots,t_{s'})^T, 
\end{equation} 
where $s'\leq s\leq m$ and $x'_1,\ldots,x'_{s'}$ are distinct.  
Since $(i, j)= (m,m+2), (m+1,m+2), (m,m+1)$, the first $s'$ rows  of $M_{s', i, j}(x'_1,\ldots, x'_{s'})$ are Vandermond determinant and so $t_i=0$ for $1\leq i\leq s'$,    
which contradicts with $t'_1(x'_1)^{m+1}+\ldots+t_{s'}(x_{s'})^{m+1}=1$, $t'_1(x'_1)^{m}+\ldots+t_{s'}(x_{s'})^{m}=1$,  or 
$t'_1(x'_1)^{m+2}+\ldots+t_{s'}(x_{s'})^{m+2}=1$, respectively.    \hfill $\Box$

\section{The girth of jumped Wenger graphs}

In graph theory, the {\it girth} of a graph is the length of a shortest cycle contained in the graph.  In  \cite{lu,SHS}, it is proved that  the Wenger graphs  have girth 8 or $6$ up to the choices of $q$ and $m$.  In \cite{clw3} authors proved that the linearized Wenger graphs have girth $8$ or $6$ depending on the choices of $q$ and $m$. 
In this section we consider the girth of the jumped Wenger graphs.

It is no harm to assume that $L_1P_1L_2\ldots P_sL_1$ is a cycle of length $2s$, where $P_h=(p^{(h)}_1, \ldots$, $p^{(h)}_{m+1} )$ and $L_h=[l^{(h)}_1,\ldots, l^{(h)}_{m+1}]$, $h=1,\ldots,s$. 

For those $L_i$, we have 
\begin{eqnarray*}
L_2-L_1 & = & [t_1, t_1f_2(p^{(1)}_1),  \ldots,  t_1f_{m+1}(p^{(1)}_1) ]\\ 
L_3-L_2 & =& [t_2, t_2f_2(p^{(2)}_1),  \ldots,  t_2f_{m+1}(p^{(2)}_1) ] \\ 
\ldots \ldots &  & \ldots \ldots \\ 
L_{s}-L_{s-1} &= &[t_{s-1}, t_{s-1}f_{2}(p^{(s-1)}_1),  \ldots,  t_{s-1}f_{m+1}(p^{(s-1)}_1) ]\\
L_1-L_{s} & = & [t_s, t_s f_2(p^{(s)}_1),  \ldots,  t_s f_{m+1}(p^{(s)}_1) ],
\end{eqnarray*} 
where $t_h=l^{(h+1)}_1-l^{(h)}_1$, for $h=1,\ldots,s-1$, $t_s=l^{(1)}_1-l^{(s)}_1$.

Therefore we have

\begin{equation} \label{eq:girth}
\left (\begin{array}{cccc}
1 & 1 & \cdots  & 1 \\
f_2(p^{(1)}_1) & f_2(p^{(2)}_1) & \cdots &  f_2(p^{(s)}_1) \\
f_3(p^{(1)}_1) & f_3(p^{(2)}_1) & \cdots &  f_3(p^{(s)}_1) \\
\multicolumn{4}{c}\dotfill \\
f_{m+1}(p^{(1)}_1) & f_{m+1}(p^{(2)}_1) & \cdots &  f_{m+1}(p^{(s)}_1) 
\end{array}
\right ) \left (\begin{array}{c}
t_1\\
t_2\\
\vdots \\
t_s
\end{array} \right )=0
\end{equation}

Again according to Fact 1 and 2, if $L_1P_1\ldots P_sL_1$ is a cycle of length $2s$, then $t_h\neq 0$ for all $h=1,\ldots, s$,  $p^{(h)}_1-p^{(h+1)}_1\neq 0$ for $h=1,\ldots,s-1$ and $p^{(s)}_1-p^{(1)}_1\neq 0$. 

In this paper we consider $(1,f_2(x),\ldots, f_{m+1}(x))=(1,x, x^2,\ldots, x^{i-1}$, $ x^{i+1},\ldots$, $x^{j-1},x^{j+1},\ldots,x^{m+2})$ and  (\ref{eq:girth}) can be written as

\begin{equation} \label{eqn:key}
M_{m+1,i, j}(p^{(1)}_1,\ldots, p^{(s)}_1)(t_1 \ldots t_s)^T=0.
\end{equation}
  
First we give an upper bound for the girth of the  jumped Wenger graphs. 

{\begin{thm} \label{thm:girth}For any integers $i, j$, $1\leq i<j\leq m+2$, the girth of the jumped Wenger graph $J_{m}(q, i, j)$  is less than or equal to $8$. 
\end{thm}} 

\noindent  \textbf{Proof.} It suffices to find a cycle of length $8$ in any graph $J_{m}(q, i, j)$ for 
$1\leq i<j\leq m+2$.  Putting $s=4$, $p^{(1)}_1=p^{(3)}_1=0$, $p^{(2)}_1=p^{(4)}_1=1$ in (\ref{eqn:key}),  the rank of $M_{m+1, i, j}(0, 1, 0, 1)$ is $2$ and $t_1+t_3=0$ and $t_2+t_4=0$.   Let $l^{(1)}_1=0$, $l^{2)}_1=1$, $l^{(3)}_1=0$, $l^{(4)}_1=1$ and $P_1=(0,\ldots,0)$. A cycle $L_1P_1L_2P_2L_3P_3L_4P_4L_1$, that is, 
\begin{eqnarray*}
P_1=(0,0,\ldots,0), & P_2=(1, 1,\ldots, 1), \\
P_3=(0,1,\ldots, 1),  & P_4=(1,0,\ldots,0), \\  
L_1=(0,0,\ldots, 0), & L_2=(1,0,\ldots,0), \\
L_3=(0,1, \ldots,1), & L_4=(1,1,\ldots,1), 
\end{eqnarray*} 
 is required. \hfill $\Box$

Next, we determine the exact girth case by case.

{\begin{thm}
If $m=1$, the girth of the jumped Wenger graph $J_m(q, i, j)$, $1\leq i<j\leq m+2$,  is  $4$ in the following cases :

(a)  $(i, j)=(1,3) \mbox{ and } char(q)\neq 2$; 

(b) $(i, j)=(1,2) \mbox{ and } 3|(q-1)$.

\noindent  For $m=1$, the girth of the jumped Wenger graph  $J_m(q, i, j)$, $1\leq i<j\leq m+2$, is $6$ in the following cases:

(a) $ (i, j)=(1,3) \mbox{ and } p=char(q)= 2$;

(b) $ (i, j)=(1,2) \mbox{ and } 3 \nmid (q-1)$;

(c) $(i, j)=(2,3)$ and $q\neq 2, 3$. 

\noindent The girth of $J_2(q, 2, 3)$ if $q=2$ or $3$ is $8$.    
\end{thm}} 

The proof is straightforward and so is omitted here.

{\begin{lem}[\cite{nied}, Lemma 6.24] \label{lem:odd}
For odd q, let $b\in \mathbb{F}_q$, $a_1, a_2\in \mathbb{F}^*_q$,  and $\eta$ be the quadratic character of $\mathbb{F}_q$.  Then the number of solutions of the equation 
 \[a_1x^2_1+a_2x^2_2=b\] 
 is $q+v(b)\eta(-a_1a_2)$, where 
$v(b)=-1$ for $b\in \mathbb{F}^*_q$ and $v(0)=q-1$. 
\end{lem}}

{\begin{lem} [\cite{nied}, Theorem 6.32 ]\label{lem:even}
 Let $\mathbb{F}_q$ be a finite field with q even and 
let $b\in \mathbb{F}_q$. 
Then for odd $n$, the number of solutions of the equation 
\[x_1x_2+x_3x_4+\ldots+x_{n-1}x_n+x^2_n=b\]
 in $\mathbb{F}^n_q$ is $q^{n-1}$. 
 For even $n$,  the number of solutions of the equation
 \[x_1x_2+x_3x_4+\ldots+x_{n-1}x_n=b\]
in $\mathbb{F}^n_q$ is $q^{n-1}+v(b)q^{(n-2)/2}$. 
 For even $n$ and $a\in \mathbb{F}_q$ with $tr_{\mathbb{F}_q/ \mathbb{F}_p}(a)=1$, 
  the
number of solutions of the equation 
\[x_1x_2+x_3x_4+\ldots+x_{n-1}x_n+x^{2}_{n-1}+ax^2_n=b\] 
in $\mathbb{F}^n_{q}$ is $q^{n-1}-v(b)q^{(n-2)/2}$.  The definition of $v(b)$ is the same as Lemma \ref{lem:odd}.  
  \end{lem}} 
  
  {\begin{lem}    \label{lem:nonzero} For distinct $a, b, c\in \mathbb{F}_q$, 
the system of equations 
\[\left \{ \begin{array}{l}
t_1+t_2+t_3=0 \\
a t_1+b t_2+c t_3=0
\end{array} \right .\]
has solution  $(t_1,t_2, t_3)$ and $t_i\neq 0$ for all $i$, $i=1, 2, 3$. 
\end{lem}}

The proof is trivial and so omitted here. 

{\begin{lem}  \label{lem:twothree} 
 There are distinct $a, b$ and $ c\in \mathbb{F}_q$ such that  $|M_{3,2,3}(a, b, c)|= 0$ if and only if  either  $q=p^l$  with $p\neq 2$ and  $q\neq 3$,  
or $q=2^{2l}$,  where $l$ is a positive integer. 
\end{lem}} 

\noindent \textbf{Proof.} 
First, if $3|\,(q-1)$, i.e.,  $q=p^{2l}$ for a positive integer $l$,  then  let $a=1, b=g^{\frac{q-1}{3}}, c=g^{\frac{2(q-1)}{3}}$, where  $g$ is  a primitive element in $\mathbb{F}_q$, and we are done. 

Next, we consider the case  $3 \nmid \,(q-1)$. If $a=0$ then $(bc^{-1})^3=1$, which implies  $b= c$.    So we can  assume that $a\neq 0$ and so we get  
\begin{equation*} 
|M_{3,2,3}(a,b,c) |=a^5\left|\begin{array}{ccc}
1 & 1 & 1 \\
0 & a^{-1}b-1 & a^{-1}c-1 \\
0 &  (a^{-1}b)^4-1 & (a^{-1}c)^4-1
\end{array}\right |.\end{equation*} 

\noindent If  $|M_{3,2,3}(a,b,c) |=0$,  we obtain  
$$\frac{(a^{-1}b)^4-1}{a^{-1}b-1}=\frac{(a^{-1}c)^4-1}{a^{-1}c-1}=-t $$ 
\noindent for some $t\in \mathbb{F}^{*}_q$ and hence $a^{-1}b$ and 
$a^{-1}c$ are different roots of the equation $x^3+x^2+x+t+1=0$.

Suppose that $x_1,x_2$ and $x_3$ are roots of $x^3+x^2+x+t+1=0$ and so we have  
\begin{equation} \label{eq:lao}
\left \{ \begin{array}{l}
x_1+x_2+x_3=-1\\
x_1x_2+x_1x_3+x_2x_3=1. \end{array} \right. 
\end{equation}

\noindent Substitute $x_3$ in Eq. (\ref{eq:lao}) by $x_3=-1-x_1-x_2$, and so we get 
\begin{equation}\label{eq:original}
x^2_1+x^2_2+x_1x_2+x_1+x_2+1=0.
\end{equation}  

\noindent Denote the set of  the solutions $(x_1, x_2)$ of Eq. (\ref{eq:original}) by $X$ and let 
\begin{equation}  
  S=\{(x_1,x_2)\in X|\, x_1\neq x_2,   x_1 \neq 0, 1, \mbox{ and } x_2\neq 0,1\}.
\end{equation}

\noindent Let us check whether $S$ is null.   Consider three cases:

 Case 1: $p\neq 2, 3$ and $q=p^{2l-1}$ where $l$ is a positive  integer.  

After multiplying $4$ on both sides of Eq. (\ref{eq:original}) and  putting $y_1=2x_1+x_2$ and  $y_2=x_2$,  
Eq. (\ref{eq:original}) becomes 
 \begin{equation} \label{eq:second}
y_1^2+3y^2_2+2y_1+2y_2+4=0. 
\end{equation}

\noindent Multiplying $3$ on both sides of Eq. (\ref{eq:second}) and putting $z_1=y_1+1$ and $z_2=3y_2+1$,    Eq. 
(\ref{eq:second}) becomes 
 \begin{equation} \label{eq:simple}
3z_1^2+z_2^2+8=0. 
\end{equation}

\noindent  By Lemma \ref{lem:odd}, the number of roots of Eq. (\ref{eq:simple}) is $q+v(-8)\eta(-3)$ and so $|X|\geq q-1$.  
 
If $x_1=x_2$, then Eq. (\ref{eq:original}) becomes 
\begin{equation}\label{eq:equal}
(3x_1+1)^2+2=0. 
\end{equation}

\noindent Thus Eq. (\ref{eq:equal}) has at most two roots.  Furthermore,  it is required that  both $x_1\neq 1$ and $x_2\neq 1$. 
If $x_1=1$,  Eq. (\ref{eq:original}) becomes $(x_2+1)^2=-2$, which has at most two roots in $\mathbb{F}_q$. The same result is for $x_2=1$. 
 However, $(x_1,x_2)=(1,1)\notin X$. 

In addition,  it is required that both $x_1$  and  $x_2$ are not zero. 
If $x_1=0$, Eq. (\ref{eq:original}) becomes $x^2_2+x_2+1=0$ and so it has two roots at most.  The result is the same if $x_2=0$. In addition, $(x_1,x_2)=(0,0)\notin X$. 

 Thus,   $|S|\geq q-1-2-4-4= q-11$.  So the result holds for any $q>11$. In addition, 
Putting   $a=1$, $b=2$ and $c=3$  if $q=5$,    $a=1$, $b=2$ and $c=-3$  if $q=7$
 and $a=-1$, $b=2$ and $c=5$ if  $q=11$, we get   $|M_{3,2,3}(a,b,c)|=0$.

Case 2. $q=2^{2l-1}$, where $l$ is a positive integer.    

Putting $y_1=x_1+1$ and $y_2=x_2+1$, 
Eq. (\ref{eq:original}) becomes 
\begin{equation}\label{eq:binary}
y_1^2+y_2^2+y_1y_2=0.
\end{equation}

\noindent According to  Lemma \ref{lem:even},  the number of solutions of Eq. (\ref{eq:binary}) in $\mathbb{F}^2_q$ is $q-v(0)q^0=1$ and so Eq. (\ref{eq:original}) has only one solution, that is, $x_1=1$ and $x_2=1$.  Therefore, there are no distinct $a$, $b$ and $c$ such that $|M_{3,2,3}(a, b, c)|= 0$. 

Case 3. $q=3^{2l-1}$, where $l$ is a positive integer.  

Multiply $4$ on both sides of Eq. (\ref{eq:original}) and  put $y_1=2x_1+x_2$ and $y_2=x_2$,  and so Eq. (\ref{eq:original}) becomes 
\begin{equation} \label{eq:triple}
(y_1+1)^2+2y_2=0.
\end{equation} 

\noindent The number of the solutions of Eq. (\ref{eq:triple}) is $q$ and so $|X|=q$.

 If $x_1=0$, Eq. (\ref{eq:original}) becomes $(2x_2+1)^2=0$ and  so $x_2=1$.  Similarly,  we have  $(1, 0)\in X$, but $(0,0)\notin X$.
If $x_1=1$,  Eq. (\ref{eq:original}) becomes $x^2_2+2x_2=0$ and  so $x_2=0$, or $x_2=1$. Thus, $(1,1)\in X$. 
Therefore, $|S|\geq q-2-1=q-3$ and so the result holds for $q>3$. 

In addition, it is east to check that there are no distinct $a$, $b$ and $c$ in $\mathbb{F}_q$ such that $|M_{3,2,3}(a, b, c)|= 0$ if $q=3$.  

So the result is established in all cases. \hfill $\Box$

{\begin{lem} \label{lem:onetwo} 
There are distinct $a, b$ and $ c\in \mathbb{F}_q$ such that  $|M_{3,1,2}(a, b, c)|= 0$ if and only if 
either  $q=p^l$  with $p\neq 2$ and  $q\neq 3$,  
or $q=2^{2l}$,  where $l$ is a positive integer.  
\end{lem}} 

\noindent \textbf{Proof.}  
The result holds for $3| (q-1)$ in the same way as the proof of Lemma \ref{lem:twothree}. 
Similarly,  assume that  $a\neq 0$ and so we obtain  
\begin{equation*}
|M_{3,1,2}(a,b,c)|=\left|\begin{array}{ccc}
1 & 1 & 1 \\
0 & (a^{-1}b)^3-1 & (a^{-1}c)^3-1 \\
0 &  (a^{-1}b)^4-1 & (a^{-1}c)^4-1
\end{array}\right |.
\end{equation*} 

\noindent We have 
\begin{equation*}
\frac{(a^{-1}b)^4-1}{(a^{-1}b)^3-1}=\frac{(a^{-1}c)^4-1}{(a^{-1}c)^3-1}=t. 
\end{equation*}

\noindent So $a^{-1}b$ and $a^{-1}c$ are the solutions of the equation 
$x^3+(1-t)x^2+(1-t)x+(1-t)=0$. Thus  $ab^{-1}$ and $ac^{-1}$ are  roots of the equation
$y^3+y^2+y+t'=0$.  Similar to the proof in Lemma \ref{lem:twothree},  the result is acquired.  \hfill $\Box$

{\begin{lem} \label{lem:onethree}
There are distinct $a,b,c\in \mathbb{F}_q$ such that $M_{3,1,3}(a^2,b^2,c^2)= 0$ if and only if $p$ is odd.   
\end{lem}}

 The proof is obvious and so is omitted  here.

By Lemma \ref{lem:nonzero}, \ref{lem:twothree}, \ref{lem:onetwo},   \ref{lem:onethree} and Theorem \ref{thm:girth}, we obtain the following result. 
{\begin{thm}
If $m=2$, the girth of the Wenger graph $J_m(q,i, j)$, $1\leq i<j\leq m+2$,  is $6$ under the following cases:

(a)  $(i, j)=(1,2)$ , $q\neq 2^{2l-1}$ and $q\neq 3$ with a positive integer $l$; 

(b)  $(i, j)=(2,3)$, $q\neq 2^{2l-1}$ and $q\neq 3$ with a positive integer $l$;

(c) $(i, j)=(1,4)$ and $q\neq 2, 3, 5$; 

(d) $(i, j)=(2,4)$ and $q\neq 2$. 

\noindent If $m=2$, the girth of $J_m(1,3)$ with $char(q)\neq 2$ is 4. 

\noindent If $m=2$, the girth of  $J_m(q, i, j)$, $1\leq i<j\leq m+2$, is $8$ under the following cases:

(a) $(i, j)=(1,2)$ and $q=2^{2l-1}$ with a positive integer $l$ or $q=3$; 

(b) $(i, j)=(2,3)$ and $q=2^{2l-1}$ with a positive integer $l$ or $q=3$;

(c) $(i, j)=(1,3)$ and $char(q)= 2$.  

(d) $(i, j)=(1,4)$ and $q= 2,3,5$ 

(e) $(i, j)=(2,4)$ and $q=2$.
\end{thm}}

\begin{lem} \label{Lem:three45} 
For $m=3,4,5$,  there are distinct $a,b,c\in \mathbb{F}_q$ such that the rank of  $M_{m+1, i, j}(a,b,c)$ is less than $3$ if and only if one of the following conditions holds:

(i)  $(i, j)=(1,4)$ and $3|\,(q-1)$. 

(ii) $(i, j)=(2,5)$ and $3|\,(q-1)$.  
\end{lem} 

The proof is clear and so is omitted here. 

\begin{thm}
For $m\geq 3$, the girth of the Wenger graph $J_m(q, i, j)$, $1\leq i<j\leq m+2$, is $8$ except the following cases:

(a)  $m=3, 4, 5$, $3|\, (q-1)$  and $(i, j)=(1,4)$. 

(b)  $m=3, 4, 5$, $3|\, (q-1)$ and $(i, j)=(2,5)$.  

(b) $m=6$, $3|(q-1)$ and $(i, j)=(2,5)$.

\noindent Furthermore, the girth of such graphs is $6$ in the above cases. 
\end{thm} 

\noindent \textbf{Proof. }  If $3|(q-1)$,   put $a=1$, $b=g^{\frac{q-1}{3}}$, and $c=g^{\frac{2(q-1)}{3}}$ where  $g$ is a primitive element in 
$\mathbb{F}_q$. So  the rank of $M_{m+1, i, j}(a,b,c)$ is less than $3$ in those cases by Lemma \ref{Lem:three45} and so the girth of those graphs is $6$ by Lemma \ref{lem:nonzero}.  For other cases,  the rank of the matrix $M_{m+1, i, j}(a,b,c)$
is $3$ and so the girth is $8$ by Theorem \ref{thm:girth}. \hfill $\Box$ 

\section{Conclusions}

In this paper we  consider the diameter and girth of  jumped Wenger graphs,  whose functions jumps at two places.  Using our similar method,   we can get similar results for  those graphs with  more jumped places.  In addition, determination of the exact diameter of any jumped Wenger graphs is our future task.

\end{document}